\documentclass[10pt,twoside]{article}
\usepackage{graphicx}
\usepackage{amsmath}
\usepackage{Latex-document}

\usepackage{amssymb}

\newtheorem{Theorem}{Theorem}
\newtheorem{Conjecture}[Theorem]{Conjecture}
\newtheorem{Definition}[Theorem]{Definition}

\newcommand{\bref}[1]{(\ref{#1})}

\newcommand{\bR}{{\mathbb{R}}}
\newcommand{\half}{\tfrac{1}{2}}

\newcommand{\scm}{\mbox{\sc{m}}}
\newcommand{\go}{\mathring{g}}
\newcommand{\PM}{\mathcal{PM}}
\newcommand{\cD}{\mathcal{D}}
\newcommand{\II}{I\!I}

\title{\bf  Mass and 3-metrics of Non-negative \vskip -2mm Scalar Curvature\vskip 6mm}

\author{Robert Bartnik\vspace*{-0.5cm}\thanks{School of Mathematics and Statistics,
University of Canberra, ACT 2601, Australia.  E-mail: robert.bartnik@canberra.edu.au}}
\date{\vspace{-8mm}}

\begin{document}
\maketitle

\thispagestyle{first} \setcounter{page}{231}

\begin{abstract}\vskip 3mm

Physicists believe, with some justification, that there should be a
correspondence between familiar properties of Newtonian gravity and
properties of solutions of the Einstein equations.  The Positive Mass
Theorem (PMT), first proved over twenty years ago
\cite{SchoenYau79b,Witten81}, is a remarkable testament to this faith.
However, fundamental mathematical questions concerning mass in general
relativity remain, associated with the definition and properties of
quasi-local mass.  Central themes are the structure of metrics with
non-negative scalar curvature, and the role played by minimal area
2-spheres (black holes).

\vskip 4.5mm


\end{abstract}

\vskip 12mm

\section{Positive Mass Theorem} \setzero

\vskip-5mm \hspace{5mm}

The Positive Mass Theorem provides a good example of \emph{``the
unreasonable effectiveness of physics in mathematics\footnote{with
apologies to Eugene Wigner \cite{Wigner60}.}''}. The need to
define mass in general relativity is motivated directly by the
physics imperative to establish a correspondence between general
relativity and classical Newtonian gravity.  Already difficulties
arise: although the vacuum Einstein equations
$Ric_{\alpha\beta}-\half Rg_{\alpha\beta}=0$ for the Lorentz
metric $g_{\alpha\beta}$ suggest (by analogy with the wave
equation, for example) that a mass (energy) which includes
contributions from the gravitational field, should be built from
the first derivatives of the field $g_{\alpha\beta}$, it is clear
that this is incompatible with coordinate invariance.

The Schwarzschild vacuum  spacetime metric, for $ r>\max(0,2\scm)$,
\begin{equation}\label{Schw}
   ds^2 = - \left(1-2\scm/r\right)\, dt^2 +
\frac{dr^2}{1-2\scm/r} +
r^2(d\vartheta^2+\sin^2\vartheta\,d\varphi^2),
\end{equation}
provides an important clue, since the parameter $\scm\in\bR$ governs the
behaviour of timelike geodesics and may be regarded as the total mass.
Note that $\scm>0$ ensures the boundary $r=2\scm$ is smooth and totally
geodesic in the hypersurfaces $t=const$.

A Riemannian 3-manifold $(M,g)$ is said to be \emph{asymptotically flat} if
$M\backslash K \simeq \bR^3\backslash B_1(0)$ for some compact $K$, and $M$
admits a  metric $\go$ which is flat outside $K$, and the metric
components $g_{ij}$ in the induced rectangular coordinates satisfy
\begin{equation}
   |g_{ij}-\go_{ij}| = O(r^{-1}), \quad |\partial_k g_{ij}|=O(r^{-2}),
\quad |\partial_k\partial_l g_{ij}|=O(r^{-3}).
\end{equation}
The total mass of $(M,g)$ is defined informally by \cite{ADM61}
\begin{equation}
    m_{ADM} = \frac{1}{16\pi}\oint_{S^2(\infty)} (\partial_i g_{ij}-\partial_j
   g_{ii})\,dS_j.
\end{equation}
If the scalar curvature $R(g) \in L^1(M)$ then $m_{ADM}$ is well-defined,
independent of the choices of rectangular coordinates and of exhaustion of
$M$ used to define $\oint_{S^2(\infty)}$ --- see
\cite{Bartnik86,Chrusciel86a,OMurchadha86} for weaker decay and smoothness
assumptions.

For simplicity, the discussion here is restricted to $C^\infty$ Riemannian
3-dimensional geometry.  This corresponds to the case of time-symmetric
initial data: $(M,g)$ is a totally geodesic spacelike hypersurface in a
Lorentzian manifold, and we can identify the local matter (equivalently,
energy) density with the scalar curvature $R(g)\ge0$.  This simplification
entails a small loss of generality: most, but not all, of the results we
describe have been extended to general asymptotically flat space-time
initial data $(M,g,K)$, where $K_{ij}$ is the second fundamental form of a
spacelike hypersurface $M$.  Some results also generalize to the closely
related Bondi mass, which measures mass and gravitational radiation flux
near null infinity, and to mass on asymptotically hyperbolic and
anti-deSitter spaces cf.~\cite{Wang01,ChruscielNagy01}, but these
involve additional complications which we will not discuss here.


The Positive Mass Theorem (PMT) in its simplest form is
\begin{Theorem}
Suppose $(M,g)$ is a complete asymptotically flat 3-manifold with
non-negative scalar
curvature $R(g)\ge0$.  Then $m_{ADM}\ge0$, and $m_{ADM}=0$ iff
$(M,g)=(\bR^3,\delta)$.
\end{Theorem}

The rigidity conclusion in the case $m_{ADM}=0$ shows that
$m_{ADM}>0$ for $(M,g)$ scalar flat (``matter-free'') but
non-flat, so $m_{ADM}$ does provide a measure of the gravitational
field.

Three distinct approaches have been successfully used to prove the PMT:
with stable minimal surfaces \cite{SchoenYau79b,SchoenYau81}; with spinors
\cite{Witten81,Lohkamp99} and the Schr\"odinger-Lichnerowicz identity
\cite{Schrodinger32,Lichnerowicz63}; and using the Geroch foliation
condition \cite{Geroch73,HuiskenIlmanen01}.  A number of other appproaches
have produced partial results: using spacetime geodesics
\cite{PenroseEt93}; a nonlinear elliptic system for a distinguished
orthonormal frame \cite{Nester88,DimakisMuller90}; and alternative
foliation conditions \cite{Jang79b,Kijowski84,Bartnik93}. The connection
between these approachs remains mysterious; the only discernable common
thread is mean curvature, and this is quite tenuous.

The application of the positive mass theorem to resolve the Yamabe
conjecture \cite{Schoen84,LeeParker87} is well known. Less well known is
the proof of the uniqueness of the Schwarzschild spacetime amongst static
metrics with smooth black hole boundary \cite{BuntingMasood87}, which we
briefly outline.

A \emph{static spacetime} is a Lorentzian 4-manifold with a
hypersurface-orthogonal timelike Killing vector.  With $V$ denoting the
length of the Killing vector, the metric $g$ on the spacelike hypersurface
satisfies the static equations
\begin{equation} \label{static}
   \begin{array}{rcl}
   Ric_g & =& V^{-1} \nabla^2 V, \\
   \Delta_g V &=& 0.
   \end{array}
\end{equation}
Smoothness implies the boundary set $\Sigma=\{V=0\}$ is totally geodesic;
analyticity of $g,V$ can be used to show the asymptotic expansions
\begin{eqnarray*}
   g_{ij} &=& (1+2m/r) \delta_{ij} + O(r^{-2}),
\\
   V &=& 1 - m/r + O(r^{-2}),
\end{eqnarray*}
as $r\to\infty$ for some constant $m\in\bR$.  The
metrics $g^{\pm} = \phi_\pm^4 g$ where $\phi_\pm = (1\pm V)/2$ both have
$R(g^\pm)=0$, and $g^+$ is asymptotically flat with vanishing ADM mass, and
$g^-$ is a (smooth) metric on a compact manifold.  Gluing two copies of
$(M,g)$ along the totally geodesic boundary $\Sigma$ and conformally
changing to $\tilde{g}=\tilde{\phi}^4g$ where $\tilde\phi=\phi_\pm$ on the
two copies of $M$, gives a complete AF manifold with $R(\tilde{g})=0$ and
vanishing mass.  The PMT shows $(\tilde{M},\tilde{g})$ is flat and it
follows without difficulty that $(M,g)$ is Schwarzschild.  This extends
previous results \cite{Israel68,Robinson77} which required the boundary to
be connected.

\section{Penrose conjecture}

\setzero\vskip-5mm \hspace{5mm}

A boundary component $\Sigma$ with mean curvature $H=0$ is called a
\emph{black hole} or \emph{horizon}, since if $(M,g)$ is a totally geodesic
hypersurface then $\Sigma$ is a trapped surface and hence, by the Penrose
singularity theorem \cite{HawkingEllis73}, lies within an event horizon and
is destined to encounter geodesic incompleteness in the predictable future.

The spatial Schwarzschild metric $ g = \frac{dr^2}{1-2\scm/r} +
r^2(d\vartheta^2+\sin\vartheta d\varphi^2) $ with $\scm<0$ shows that the
completeness condition in the PMT is important, but it can be weakened to
allow horizon boundary components of $M$.  This follows immediately from
the minimal surface argument \cite{SchoenYau79b}; or by an extension to the
Witten argument \cite{GHHP83}, imposing one of the boundary conditions
\begin{equation}
   \psi = \pm \epsilon \psi \textrm{   on  }\Sigma=\partial M,
\end{equation}
on the spinor field $\psi$, where $\epsilon=\gamma^n \gamma^0$ satisfies
$\epsilon^2=1$.  An interesting extension is obtained by imposing the spectral boundary condition
\begin{equation}\label{P+psi}
   P_+ \psi = 0  \textrm{   on  }\Sigma
\end{equation}
where $P_+$ is the projection onto the subspace of positive eigenspinors of
the induced Dirac operator $\cD_\Sigma$.  Using the remarkable Hijazi-B\"ar
estimate \cite{Hijazi86,Bar92}
\begin{equation}\label{Hijazi}
   |\lambda| \ge \sqrt{{4\pi}/{|\Sigma|}},
\end{equation}
for the eigenvalues of $\cD_\Sigma$ when $\Sigma\simeq S^2$, Herzlich
showed \cite{Herzlich97a}
\begin{Theorem}
  If $(M,g)$ is asymptotically flat with $R(g)\ge0$ and boundary
  $\Sigma\simeq S^2$ with mean curvature satisfying
  \begin{equation}
     H_\Sigma \le 2/r
  \end{equation}
  where $r= \sqrt{|\Sigma|/4\pi}$, then $m_{ADM}\ge 0$, with equality
  iff $(M,g)= (\bR^3\backslash B(r),\delta)$.
\end{Theorem}

The proof starts with the Riemannian form of the
Schr\"odinger-Lichnerowicz-Witten identity
\cite{Schrodinger32,Lichnerowicz63,Witten81}
\begin{equation}\label{SLidentity}
    \int_M(|\nabla\psi|^2+\tfrac{1}{4}R(g)|\psi|^2-|\cD\psi|^2)\,dv_M
    = 4\pi|\psi_\infty|^2 m_{ADM} +  \oint_\Sigma \mu(\psi),
\end{equation}
where $\mu(\psi)$ is the Nester-Witten form \cite{Nester83}
\begin{eqnarray}
   \mu(\psi)
\label{NWform}
 &=&    \langle\psi,(\cD_\Sigma+\half H_\Sigma)\psi\rangle\, dv_\Sigma.
\end{eqnarray}
The boundary condition $P^+\psi|_\Sigma=0$ is
elliptic and it can be shown \cite{BartnikChrusciel02} there is a spinor on
$M$ satisfying $\cD\psi=0$ with boundary conditions
$\psi\to\psi_\infty\ne0$ as $r\to\infty$ and \bref{P+psi} on $\Sigma$.  It
follows from \bref{Hijazi} and \bref{P+psi} that
 $   \langle\psi,(\cD_\Sigma+\half H_\Sigma)\psi\rangle
\le (\half H_\Sigma - |\lambda_1^-|)|\psi|^2 \le 0
$
and the result follows.

Observe that in each case, equality leads to flat $\bR^3$.
An elegant physical argument lead Penrose to conjecture an analogous
inequality, but which distinguishes the Schwarzschild metric instead
\cite{Penrose73}, see also \cite{Gibbons84}.
\begin{Conjecture}[Penrose]
  If $(M,g)$ satisfies the conditions of the PMT, except that $\partial
  M=\Sigma$ is compact with vanishing mean curvature and such that $\Sigma$
  is the ``outermost'' closed minimal surface in $M$, then
\begin{equation}\label{Penrose}
  m_{ADM} \ge
  \sqrt{{|\Sigma|}/{16\pi}},
\end{equation}
  with equality only for the  Schwarzschild metric.
\end{Conjecture}

A closed minimal surface is said to be an \emph{outermost horizon} or
\emph{outer-minimizing horizon} if $M$ contains no least area surfaces
homologous to $\Sigma$ in the asymptotic region exterior to $\Sigma$.  The
outermost condition is essential, since examples of non-negative scalar
curvature manifolds can be constructed by forming the connected sum of $M$
and large spheres by arbitrarily small and large necks.

The Penrose conjecture has been established by Huisken and Ilmanen
\cite{HuiskenIlmanen97,HuiskenIlmanen01} using a variational level set
formulation of the inverse mean curvature flow \cite{Geroch73}, and by Bray
\cite{Bray01} by a very interesting conformal deformation argument. Bray's
proof is more general since it takes into account contributions from all
the connected components of the boundary.

\section{Quasi-local mass} \setzero

\vskip-5mm \hspace{5mm}

Thus it is natural to consider $\sqrt{|\Sigma|/16\pi}$ as the mass of a black
hole (minimal surface) $\Sigma$.  More generally, the correspondence with
Newtonian gravity suggests that any bounded region $(\Omega,g)$ should have
a \emph{quasi-local} mass, which measures both the matter density
(represented in this case by the scalar curvature $R(g)\ge0$), and some
contribution from the gravitational field.
The rather satisfactory positivity properties of the total mass, as
established by the PMT, motivate the  properties we might expect
such a geometric mass to possess
\cite{Eardley79,ChristodoulouYau88,Bartnik97d}.
\begin{enumerate}
\item {\textbf{(non-negativity)}} \ \
         $m_{QL}(\Omega)\ge0$;
\item {\textbf{(rigidity/strict positivity)}} \ \
         $m_{QL}(\Omega)=0$ if and only if $(\Omega,g)$ is flat;
\item {\textbf{(monotonicity)}}\ \
         $m_{QL}(\Omega_1)\le m_{QL}(\Omega_2)$ whenever
$\Omega_1\subset \Omega_2$, where it is understood that the inclusion is a
metric isometry;
\item {\textbf{(spherical mass)}}\ \
         $m_{QL}$ should agree with the spherical mass, for
spherically symmetric regions;
\item {\textbf{(ADM limit)}}\ \
         $m_{QL}$ should be asymptotic to the ADM mass;
\item {\textbf{(black hole limit)}}\ \
         $m_{QL}$ should agree with the black hole mass \bref{Penrose}.
\end{enumerate}
Many candidates have been proposed for quasi-local mass (see
for example \cite{Bergqvist93} for a comparison of some definitions), the
most significant being that of Hawking
\cite{Hawking68},
\begin{equation}
   m_H(\Sigma) =
   \sqrt{\frac{|\Sigma|}{16\pi}}
   \left(1-\frac{1}{16\pi}\oint_{\Sigma}H^2\right)
\end{equation}
where $\Sigma=\partial \Omega$. This equals $\scm$ for
standard spheres in Schwarschild.  Although $m_H\le0$ for surfaces in
$\bR^3$, it was shown in \cite{ChristodoulouYau88} that $m_H(\Sigma)\ge0$
for a stable constant mean curvature 2-sphere $\Sigma$ in a 3-manifold of
non-negative scalar curvature.  Thus for such ``round'' spheres, $m_H$ is
nonegative, and the black hole limit condition is trivially satisfied.
However the remaining properties, in particular rigidity and monotonicity,
are rather problematic.  Although the twistorially-defined Penrose
quasi-local mass \cite{Penrose82} is well-behaved in special cases
\cite{Tod83}, it is defined unambiguously only for surfaces arising from
embedding into a conformally flat spacetime, and even then numerical
experiments \cite{BernsteinTod94} strongly suggest that monotonicity is
violated.

In fact, of the various proposals for $m_{QL}$, only the definitions of
\cite{ChristodoulouYau88,Bartnik88,DouganMason91} are known
to satisfy positivity.  Dougan and Mason \cite{DouganMason91} show that
the integral $\oint_\Sigma\mu(\psi)$ of the Nester-Witten 2-form
\bref{NWform} is positive for spinor fields $\psi$ on $\Sigma$ which
satisfy a certain elliptic system on $\Sigma$.  However, Bergqvist
\cite{Bergqvist92a} shows that positivity holds under much weaker
conditions on $\psi$, and there are many variant definitions with similar
properties.  It would be useful to understand these DM-style definitions
better, and in particular whether any satisfy monotonicity.

Monotonicity and ADM-compatibility  imply $m_{QL}(\Omega)\le
m_{ADM}(M,g)$ for any region $\Omega$ embedded isometrically in an
$(M,g)$ satisfying (as always) the PMT
conditions.  This motivates the following definition
\cite{Bartnik89a,HuiskenIlmanen01}
\begin{Definition} \label{qlmass}
  Let $\PM$ denote the set of all asymptotically flat 3-manifolds $(M,g)$
  of non-negative scalar curvature, with boundary which if non-empty,
  consists of compact outermost horizons, and such that $(M,g)$ has no
  other horizons.  For any bounded open connected region
  $(\Omega,g)$, let $\PM(\Omega)$ be the set of $(M,g)\in\PM$ such that
  $\Omega$ embeds isometrically into $M$, and define
  \begin{equation}
   m_{QL}(\Omega) = \inf \{ m_{ADM}(M,g) : (M,g)\in \PM(\Omega)\}.
  \end{equation}
   We say that $M$ satisfying these conditions  is an \emph{admissible
  extension} of $\Omega$.
\end{Definition}

The horizon condition serves to exclude examples which hide
$\Omega$ inside an arbitrarily small neck, which would force the
infimum to zero. This is a refinement \cite{HuiskenIlmanen01} of
the original definition \cite{Bartnik89a}, which prohibited
horizons altogether.

Clearly $m_{QL}(\Omega)$ is well-defined and finite, once the region
$\Omega$ admits just one admissible extension.  The PMT with horizon
boundary implies non-negativity, and monotonicity follows directly.  Strict
positivity of $m_{QL}$ was established in \cite{HuiskenIlmanen01}, with the
slightly weaker rigidity conclusion that if $m_{QL}(\Omega)=0$ then
$\Omega$ is locally flat.  Agreement with the spherical mass, and the ADM
limit condition, follows also from \cite{HuiskenIlmanen01}.  Bray's results
imply that $m_{QL}(\Omega)$ agrees with the black hole mass in the limit as
$\Omega$ shrinks down to a black hole.  In addition, $m_{QL}(\Omega) \le
m_{ADM}(M)$ for any admissible extension $M$, so $m_{QL}$ is the optimal
quasi-local mass definition with respect to this condition.

The optimal form of the horizon condition remains conjectural.  Bray has
suggested an alternative condition, that $\Omega$ be a ``strictly
minimizing hull'' \cite{HuiskenIlmanen01} in $M$, so $\Sigma
=\partial\Omega$ has the least area amongst all enclosing surfaces in the
exterior.  In this case we say $\Sigma$ is \emph{outer minimizing}, and
denote by $\tilde{m}_{QL}(\Omega)$ the quasilocal mass function defined by
restricting admissible extensions to those $M$ in which $\Sigma$ is outer
minimizing.  For this modified definition the Penrose inequality
\cite{HuiskenIlmanen01,Bray01} applies to show that if $\partial\Omega$
embeds into the Schwarzschild 3-manifold with the same induced metric and
mean curvature (cf. \bref{2ndvar}, \bref{staticbndry}) and encloses the
horizon, then $m_{QL}(\Omega) = \scm$.  It is not clear how to establish
this natural result for the unmodified definition $m_{QL}(\Omega)$.

\section{Static metrics} \setzero

\vskip-5mm \hspace{5mm}

Although in many respects the definition of $m_{QL}$ is quite satisfactory,
it is not constructive, and thus it is important to determine computational
methods.  The key is the following \cite{Bartnik89a}
\begin{Conjecture}\label{staticconjecture}
   The infimum in $m_{QL}$ is realised by a 3-metric agreeing with
   $\Omega$ in the interior, \emph{static} \bref{static} in the exterior
   region, and such that the metric is Lipschitz-continuous across the
   matching surface $\Sigma$, and the mean curvatures of the two sides
   agree along $\Sigma$.
\end{Conjecture}

A similar conjecture for the space-time generalisation of the quasi-local
mass, asserts that the exterior metric is \emph{stationary}, ie.~admits a
timelike Killing field \cite{Bartnik89a,Bartnik97d}.

As motivation for this conjecture, note first that if $R(g)>0$ in some
region, then a conformal factor $\phi$ can be found such that $\phi^4g$ has
less mass and $R(\phi^4g)\ge0$.  Thus a mass-minimizing metric for
\bref{qlmass}, if such a metric exists, must have vanishing scalar
curvature. Now if the linearization $ DR(g) h = \delta_g\delta_g h - \Delta
tr_gh - Ric\cdot h $ is surjective then $g$ admits a variation which
produces positive scalar curvature.  The formal obstruction to surjectivity
is non-trivial $\ker DR(g)^*$, which leads to the static metric equations
\bref{static}.  Corvino \cite{Corvino00} shows that if $\ker DR(g)^*$ is
trivial in $U\subset M$ then there are compactly supported metric
variations in $U$ which increase the scalar curvature. This gives
\begin{Theorem}
If $(M,g)$ realizes the infimum in Definition \ref{qlmass}, then there is a
$V\in C^\infty(M\backslash\Omega)$ such that $g,V$ satisfy the
static metric equations \bref{static} in $M\backslash\Omega$.
\end{Theorem}

This suggests a computational algorithm for determining $m_{QL}(\Omega)$:
find an asymptotically flat static metric with boundary geometry matching
that of $\partial\Omega$. To determine the appropriate boundary conditions,
recall the second variation formula for the area of the leaves of a
foliation labelled by $r$:
\begin{equation} \label{2ndvar}
  R(g) = 2D_nH - |\II|^2 - H^2 + 2K - 2 \lambda^{-1}\Delta_r\lambda
\end{equation}
where $\II,H,K$ are respectively the second fundamental form, mean
curvature and Gauss curvature of the leaves, $\lambda$ is the lapse
function, $n=\lambda^{-1}\partial_r$ is the normal vector and $\Delta_r$ is
the Laplacian on the leaves.  Our conventions give $H=-D_n(\log\sqrt{\det
g_r})$ where $g_r$ is the volume element of the leaves.
This shows that $R(g)$ will be defined distributionally across a matching
surface as a bounded function if
\begin{equation} \label{staticbndry}
   \begin{array}{rcl}
   g|_{T\partial\Omega}&=& g|_{T\Sigma}, \\[3pt]
   H_{\partial\Omega}&=& H_{\Sigma} .
   \end{array}
\end{equation}
\begin{Conjecture}\label{staticexistence}
$(\Omega,g)$ determines a unique static asymptotically flat manifold
$(S,g)$ with boundary $\Sigma\simeq\partial\Omega$ satisfying
\bref{staticbndry}.
\end{Conjecture}

If true, this would give a prime candidate for the minimal mass extension.
It is known (Pengzi Miao, private communication) that the boundary
conditions \bref{staticbndry} are elliptic for \bref{static}.

It is tempting to conjecture that mass-minimizing sequences for $m_{QL}$
should converge to a static metric.  For example, \cite[Theorem
5.2]{Bartnik86} shows that a sequence of metrics $g_k$, close in the
weighted Sobolev space $W^{2,q}_{-\tau}$, $q>3,\tau>1/2$, to the flat
metric $\delta$ on $\bR^3$ and such that $m_{ADM}(g_k)\to0$, converges
strongly to $\delta$ in $W^{1,2}$.  Similar results, under rather different
size conditions, are given in \cite{FinsterKath01}, and a discussion of the
general ``weak compactness'' conjecture may be found in
\cite{HuiskenIlmanen01}.

\section{Estimating quasi-local mass}

\setzero\vskip-5mm \hspace{5mm}

To estimate $m_{QL}$ from above, it suffices to construct admissible
extensions --- metrics with non-negative scalar curvature and satisfying
\bref{staticbndry}.  These boundary conditions  exclude the usual
conformal method.  Instead, metrics in \emph{quasi-spherical} form
\cite{Bartnik93}
\begin{equation}\label{QSmetric}
   g=u^2\,dr^2 + (r\,d\vartheta+\beta^1dr)^2 +
   (r\sin\vartheta\,d\varphi+\beta^2dr)^2
\end{equation}
satisfy a parabolic equation for $u$ on $S^2$ evolving in the radial
direction, when $R(g)=0$, with $\beta^1,\beta^2$ freely specifiable.  Since
the metric 2-spheres $S^2_r$ have mean curvature $H_r =
(2-\mathrm{div}_{S^2}\beta)/ur >0$, \bref{QSmetric} provides admissible
extensions for $\partial\Omega = S^2_r$ with mean curvature $H>0$. The
underlying parabolic equation derives from \bref{2ndvar}, and has been
generalized to non-spherical foliations in \cite{SmithWeinstein00}.  As an
application, choosing $\beta=0$ we can show
\begin{Theorem}
Suppose $\partial\Omega =S^2_r$ metrically, with $H\ge0$. Then
\begin{equation}
   m_{QL}(\Omega) \le \half r ( 1-\tfrac{1}{4}r^2\min_{\partial\Omega}H^2).
\end{equation}
\end{Theorem}

This bound is sharp when $\Omega$ is a flat ball or a Schwarzschild
horizon.

Finding lower bounds for $m_{QL}(\Omega)$ is more difficult.  Bray's
definition of \emph{inner mass} \cite[p243]{Bray01} gives a lower bound,
but for $\tilde{m}_{QL}(\Omega)$.  The difficulty here as above lies in
showing that a horizon inside $\Omega$ remains outermost when the inner
region is glued to a general exterior region $M_{\mathrm{ext}}\subset
M\in\PM(\Omega)$.  This follows easily when $\Sigma=\partial\Omega$ is
outer-minimizing in $M_{\mathrm{ext}}$, as guaranteed by the definition for
$\tilde{m}_{QL}(\Omega)$.

On physical grounds one expects that if ``too much'' matter is compressed
into region which is ``too small'', then a black hole must be present.  The
geometric challenge lies in making this heuristic statement precise, and
the only result in this direction has been \cite{SchoenYau83}, which gives
quantitative measures which guarantee the existence of a black hole.  An
observation by Walter Simon (private communication) is thus very
interesting: if $m_{QL}(\Omega)=1$ (say) and $\Omega$ embeds isometrically
into a complete asymptotically flat manifold $M$ without boundary and with
non-negative scalar curvature, and such that $m_{ADM}(M)<1$, then $M$ must
have a horizon.  This reinforces the importance of finding good lower
bounds for $m_{QL}$, since the existence of a horizon in a similar
situation with $\tilde{m}_{QL}$ does not follow.

\label{lastpage}

\end{document}